\documentclass[12pt]{amsart}
\usepackage{amssymb, amstext, amscd, amsmath}
\theoremstyle{plain}
\newtheorem{thm}{Theorem}[section]
\newtheorem{cor}[thm]{Corollary}
\newtheorem{prop}[thm]{Proposition}
\newtheorem{lem}[thm]{Lemma}
\theoremstyle{definition}
\newtheorem{rem}[thm]{Remark}
\newtheorem{defn}[thm]{Definition}

\newtheorem{eg}[thm]{Example}
\newcommand{\bsl}{\setminus}

\newcommand{\cstar}{\mathrm{C}^*}

\newcommand{\ol}{\overline}
\newcommand{\td}{\widetilde}
\newcommand{\wh}{\widehat}

\newcommand{\scp}{{\mathbb{Z}_+^2 \times_{\ga}}}
\newcommand{\cp}{{\mathbb{Z}^2 \times_{\ga}}}

\newcommand{\bbC}{{\mathbb{C}}}
\newcommand{\bbD}{{\mathbb{D}}}

\newcommand{\bbT}{{\mathbb{T}}}
\newcommand{\bbZ}{{\mathbb{Z}}}

  \newcommand{\A}{{\mathcal{A}}}
  \newcommand{\B}{{\mathcal{B}}}
  
  \newcommand{\D}{{\mathcal{D}}}

\renewcommand{\H}{{\mathcal{H}}}
  \newcommand{\I}{{\mathcal{I}}}
  \newcommand{\J}{{\mathcal{J}}}
  \newcommand{\K}{{\mathcal{K}}}
\renewcommand{\L}{{\mathcal{L}}}
  \newcommand{\M}{{\mathcal{M}}}
  
\renewcommand{\O}{{\mathcal{O}}}

\newcommand{\iAk}{{\mathit{A}_k}}

\newcommand{\iAd}{{\mathit{A}(\mathbb{D})}}
\newcommand{\iAdd}{{\mathit{A}(\mathbb{D}^2)}}
\newcommand{\iAkd}{{\mathit{A}_k(\mathbb{D})}}
\newcommand{\iAkdd}{{\mathit{A}_k(\mathbb{D}^2)}}
\newcommand{\iN}{{\mathit{N}}}
\newcommand{\il}{{\mathit{l}}}

\newcommand{\gl}{{\lambda}}
\newcommand{\go}{{\omega}}
\newcommand{\ga}{{\alpha}}
\newcommand{\gd}{{\delta}}
\newcommand{\gs}{{\sigma}}
\newcommand{\gp}{{\psi}}
\newcommand{\gt}{{\tau}}

\newcommand{\Ad}{\operatorname{Ad}}

\newcommand{\ba}{\backslash}
\begin{document}

%%%%%%%%%%%%%%%%%%%%%%%%%%%%%%%%%%%%%%%%%%%%%%%%%%%%%%%%%%%%
\title[Semicrossed Products]{Semicrossed products generated by two commuting automorphisms}

\author[M. Alaimia]{Mohamed Ridha Alaimia}
\thanks{}
\address{Department of Mathematical Sciences\\
King Faud University of Petroleum and Minerals\\
Dhahran 31261, Saudi Arabia}
\email{alaimia@kfupm.edu.sa}

\author[J. Peters]{Justin R. Peters}
\address{Department of Mathematics\\
Iowa State University\\
Ames, IA
50011}
\email{peters@iastate.edu}

\begin{abstract}
In this paper, we study the semicrossed product of a finite dimensional
$\cstar$-algebra for two types of $\bbZ_{+}^2 $-actions, and identify them
with matrix algebras of analytic functions in two variables. We look at the connections with
semicrossed by $\bbZ_{+}$-actions.

\end{abstract}

\subjclass{47L80}
\date{}
\maketitle
%%%%%%%%%%%%%%%%%%%%%%%%%%%%%%%%%%%%%%%%%%%%%%%%%%%%%%%%%%%%

\section{Introduction} \label{s:intro}
The study of semicrossed products with respect to $\bbZ_{+}$-actions was begun in
\cite{arveson}, \cite{peters}, and higher dimensional actions have also been
considered (e.g., \cite{lingmuhly}, \cite {power}). In this note we continue the
program of \cite{dealba}, \cite{alaimia1}, and \cite{alaimia2}, which studied
$\bbZ_{+}$-actions on finite dimensional $\cstar$ algebras, to the case of
$\bbZ_{+}^2$-actions.

For a transitive action $\gs$ on a finite set $X$, card$(X) = k$, the crossed product $\bbZ \times_{\gs} C(X)$
is $M_k \otimes C(\bbT)$. The semicrossed product, which is the closed subalgebra corresponding to nonnegative
powers of the automorphism, is not $M_k \otimes \iAd$ (where $\iAd$ is the disk algebra), but rather a proper subalgebra
(cf section \ref{s:matrixfunctions}). For a transitive $\bbZ^2$-action on $X$, the crossed product is
$M_k \otimes C(\bbT^2)$, and the semicrossed product is a matrix algebra of bianalytic functions,
which is a proper subalgebra of
$M_k \otimes \iAdd$, (where $\iAdd$ is the bidisk algebra), card$(X) = k$.

However, different transitive actions, though having the same $\cstar$-crossed product, may give rise to
nonisomorphic semicrossed products.  In fact, we exhibit two transitive $\bbZ^2$-actions $\gs,\ \gs'$ on a finite
set $X$ for which the semicrossed products are not isomorphic.

In Proposition \ref{p:z2} we show that if $\gs = (\gs_1, \gs_2)$ is a $\bbZ^2$-action
on $X,\ \text{card}(X) = k$, then the semicrossed product with respect to the action
on $C(X)$ is the norm closure of the algebra generated by
\[ zP_1, \quad wP_2, \quad  D \in \D, \]
where $P_i$ is the $k \times k$ permutation matrix associated with $\gs_i, \ i = 1, 2$, and $\D$ is the algebra of
diagonal matrices in $M_k$. The closure is in $M_k \otimes \iAdd $.

Ideally, a classification theorem for these $\bbZ^2$-actions on $X$ would establish
necessary and sufficient conditions for two semicrossed products to be (completely)
isometrically isomorphic, in terms of some invariants of the permutation matrices.
Such a theorem has not been achieved here.  Our results, which enable us to
distinguish the semicrossed products coming from certain (non-conjugate)
$\bbZ^2$-actions, are obtained by computing the codimensions of maximal ideals of the
algebra.

One of the main theorems for semicrossed products of the form
\mbox{$\bbZ_{+} \times_{\gs} C(X)$} (for $X$ a compact metric
space) is that two such semicrossed products are isomorphic if and only if the actions are conjugate.  For
$\bbZ^2$-actions, the question is open.  The work here lends support to the conjecture that the analogue of the
theorem for $\bbZ$-actions may hold for $\bbZ^2$-actions as well.

\section{Some matrix function algebras} \label{s:matrixfunctions}

We will let $\bbD$ be the open unit disk, and $\bbD^2 = \bbD \times \bbD,$ the bidisk.  $\iAd $ will denote the disk
algebra, the subalgebra of $C(\overline{\bbD})$ of functions which are analytic on
the interior of $\bbD,$ and $\iAdd$ the
bidisk algebra; i.e., the subalgebra of $C(\overline{\bbD}^2)$ of bianalytic functions. We recall the matrix algebras of
analytic functions studied in \cite{dealba}:

\[ \B_k = \{ (f_{ij})^{(k-1}_{i,j=0} \ : f_{ij} \in \iAd,\ f_{ij} \sim \sum_{n=0}^{\infty} a_n^{(ij)}z^{l(i,j) + nk} \} \]

where $l(i,j) \in \{0, \dots, k-1 \}$ and $l(i,j) \equiv i - j [\textrm{mod} k].$  Note we are labeling the matrix entries
from $0$ to $k-1$ rather than from $1$ to $k.$

For $k \in \bbZ^+,$
let $\iAk$ denote the subalgebra of $\iAd $ whose nonzero Fourier coefficients are multiples of $k.$
Then the algebra $\B_k$ takes the form

\begin{equation} \label{E:Ak}
\begin{pmatrix}
\iAk & z^{k-1}\iAk & \cdots & z\iAk \\
z\iAk & \iAk & \cdots & z^2\iAk \\
\vdots & \vdots & \ddots & \vdots \\
z^{k-1}\iAk & z^{k-2}\iAk & \cdots & \iAk
\end{pmatrix}
\end{equation}

The algebras $\B_k$ were first studied in \cite{dealba} and later in \cite{alaimia1}, \cite{alaimia2},
arise as semi-crossed products of actions of $\bbZ^+$ on
$\bbC^k .$

 Next we introduce a two variable analog of this algebra.

\[ \B_{k,2} = \{ (f_{ij})_{i,j = 0}^{k-1}: \ f_{ij} \in \iAdd ,\ f_{ij} \sim \sum_{\substack{m,n \geq 0, \\
 m+n \equiv l(i,j) [\textrm{mod} k]}} ^{\infty} a_{mn}^{(i,j)} z^m w^n \} \]

Another algebra which will play an auxiliary role is

\[ \iAkdd := \{ f \in \iAdd:\ f \sim \sum_{\substack{m, n \geq 0,  \\ m+n = 0 [\text{mod} k]}}^\infty
 a_{mn}z^mw^n   \} .\]

Using Fejer's theorem in two variables (\cite[vol 2, p. 304]{zygmund}) and the fact that Cesaro means of a funtion in
$\iAkdd$ are polynomials in the same algebra, we see that the polynomials in $\iAkdd$ are dense in $\iAkdd.$
Similarly, the subalgebra of $\B_{k,2}$ whose matrix entries are polynomials is dense in $\B_{k,2}.$

\section{Maximal Ideals} \label{s:maxid}

One of the main tools in studying and distinguishing various matrix function algebras
is the space of maximal ideals. First we recall (\cite{alaimia1}, \cite{alaimia2}) the
maximal ideal structure of $\B_{k}$. Let $F = (f_{ij}) \in \B_{k}$ and $0 \leq i_0
\leq k-1$. Set $\gp_0^{i_0} : \B_{k} \to \bbC,\ \gp_0^{i_0}(F) = f_{i_0i_0}(0),$ and
for $ \gl \in \ol{\bbD}\bsl \{0\},$ set $\gp_{\gl}(F) = F(\gl)$, that is,
$\gp_{\gl}(F)$ is the $k \times k$ matrix $(f_{ij}(\gl))$. Then $J_0^{i_0} :=
\text{ker}(\gp_{0}^{i_0})$, and $J_{\gl} := \text{ker}(\gp_{\gl})$, $\gl \in \ol{\bbD}
\bsl \{0\}$ are all the maximal ideals of $\B_{k}$.

To study the maximal ideals of $\B_{k,2}$ we first look at $\iAkdd.$

\begin{lem} \label{l:iAkddmaxid}
$\iAkdd$ is a Banach subalgebra of $\iAdd$ and its maximal ideals are of the form

\[ \iN_{\gl, \mu} = \{ f \in \iAkdd: \ f(\gl, \mu) = 0 \} \text{ for } (\gl, \mu) \in \overline{\bbD}^2 . \]
\end{lem}

\begin{proof} It is routine to check that $\iAkdd$ is a Banach subalgebra. For the maximal ideals, let $\chi$ be a
multiplicative linear functional on $\iAkdd.$ Let $\gl$ be a $k$-th root of $\chi(z^k)$ and $\mu$ a $k$-th root of
$\chi(w^k).$ Note that $\chi(zw^{k-1})$ has the form $\go \gl \mu^{k-1},$ where $\go$ is a $k$-th root of unity.  Indeed,
$(\chi(zw^{k-1}))^k = \chi(z^k)(\chi(w^k))^{k-1} = (\gl \mu^{k-1})^k.$  Thus, replacing $\gl$ by $\go \gl$ we may assume that
$\chi(zw^{k-1}) = \gl \mu^{k-1}.$  Then

\[ \chi(z^2w^{k-2}) = (\chi(zw^{k-1}))^2 \chi(w^k)^{-1} = \gl^2 \mu^{k-2} .\]

Continuing in this way, we have that $\chi(z^rw^{k-r}) = \gl^r \mu^{k-r}, \ 0 \leq r \leq k.$  Hence if $p$ is any
polynomial in $\iAkdd,\ \chi(p) = p(\gl, \mu).$  By continuity this holds for all functions in $\iAkdd.$
\end{proof}

\begin{rem}
For $k > 1,$ the correspondence between maximal ideals and points in the bidisk is not bijective.  Indeed, for
$(\gl, \mu) \in \overline{\bbD}^2$ and $\go$ a $k$-th root of unity, the maps $f \in \iAkdd \mapsto f(\gl, \mu)$ and
$f \in \iAkdd \mapsto f(\go\gl, \go\mu)$ are identical.  Only in the case $k=1,$ the bidisk algebra, is the
correspondence bijective.
\end{rem}

Next we identify the maximal ideals of $\B_{k,2}.$  Let $F = (f_{ij})_{i,j = 0}^{k-1} \in \B_{k,2}.$ Observe that the
map $F \mapsto F(0,0)$ maps $F$ to a $k\times k$ matrix which is zero except along the diagonal. Thus if we fix
$\ell,\ 0 \leq \ell \leq k-1,$ the map

\[ \psi_0^{\ell}:  \B_{k,2} \to \bbC, \quad \psi_0^{\ell}(F) = f_{\ell \ell}(0, 0), \quad F = (f_{ij})_{i, = 0}^{k-1} \]

is a (nonzero) multiplicative linear functional, and hence $ker(\psi_0^{\ell})$ is a maximal ideal.  We will call this
ideal ``of zero type'' and denote it by $\J_{(0,0)}^{\ell}.$

A second type of mapping is parameterized by $(\gl, \mu) \in \overline{\bbD}^2\ba \{ (0, 0) \} $ and takes the form

\[ \psi_{\gl, \mu} : \B_{k,2} \to \M_k, \quad \psi_{\gl, \mu}(F) = F(\gl, \mu) .\]

By this we mean of course that each of the component functions of $F$ is evaluated at the point $(\gl, \mu).$  It is
routine to check that this mapping is surjective and multiplicative. Thus $ker(\psi_{\gl, \mu})$ provides us with
another type of maximal ideal, which we denote $\J_{(\gl, \mu)}.$

\begin{prop} \label{p:bk2maxid}
Let $\J$ be a maximal ideal of $\B_{k,2}.$  Then either
\begin{align*}
\J &= \J_{(\gl, \mu)}, \text{ for some } (\gl, \mu) \neq (0, 0) \text{ in the closed bidisk, or} \\
\J &= \J_{(0,0)}^{\ell} \text{ for some } \ell \in \{0, \dots, k-1\}.
\end{align*}
\end{prop}

\begin{proof}
Let $e_{ij},\ 0 \leq i,j \leq k-1$ be matrix units for $\M_k.$ That is, $e_{ij}$ is the matrix which is $1$ in the
$(i, j)$ position, and zero elsewhere. Since $\J$ is proper, it cannot contain all the $e_{ii},$ for then it would
contain the identity. We will assume $\J$ does not contain $e_{00}\otimes 1;$ the other cases are similar.

For any $F = (f_{ij}) \in \J$ we have

\[ (e_{00}\otimes 1)F(e_{j0}\otimes z^j) = e_{00}\otimes z^jf_{0j}, \]

so $e_{00}\otimes z^jf_{0j} \in \B_{k,2}\J\B_{k,2} = \J.$\\
Similarly, $e_{00}\otimes w^jf_{0j} \in \J.$  Also,

\[ e_{00}\otimes z^{k+j-i}f_{ij} = (e_{0i}\otimes z^{k-i})F(e_{j0}\otimes z^j) \in \J \text{ for } \in \{0, \dots, k-1 \} .\]

Similarly, $e_{00}\otimes w^{k+j-i}f_{ij} \in \J.$

Consider the closed ideal $\I$ in $\iAkdd$ generated by 
\[z^jf_{0j}, w^jf_{0j}, z^{k+j-i}f_{ij} \text{ and } w^{k+j-i}f_{ij} \text{ for }0 \leq i, j \leq k-1, \]
and where the $f_{ij}$ are the entries of an  $F =(f_{ij})$
as $F$ ranges through $ \J.$  Note that
$e_{00} \otimes \I \subset \J.$  Since by assumption $\J$ does not contain $e_{00} \otimes 1,$ it follows that
$1 \notin \I,$ and $\I$ is proper in $\iAkdd.$ By Lemma \ref{l:iAkddmaxid} , $\I \subset
\{ f \in \iAkdd:\ f(\gl, \mu) = 0 \}$ for some $(\gl, \mu) \in \overline{\bbD}^2.$

If $(\gl, \mu) = (0,0),$ then $\J \subset \J^0_{00}.$
Indeed, elements of $\J^0_{00}$ have the same components as elements of
$\B_{k,2},$ except for the $(0,0)$ component.  As $\J$ is maximal, this forces $\J = \J^0_{00}.$

If $\gl \neq 0$ and $F = (f_{ij}) \in \J,$
then $ \gl^jf_{0j}(\gl, \mu) = 0$ and $\gl^{k+j-i}f_{ij}(\gl, \mu) = 0,$ which imply $F(\gl, \mu) = 0.$
Similarly, if $\mu \neq 0,$ we conclude $F(\gl, \mu) = 0.$   But then
$\J \subset \J_{\gl,\mu},$ and so by maximality $ \J = \J_{\gl, \mu}.$
\end{proof}

\section{Semicrossed Products by $\bbZ^2_+$-actions} \label{s:scp}
Let $\A$ be a $\cstar$-algebra, and $Aut(\A)$ the group of star automorphisms of $\A.$ An action $\ga$ of $\bbZ^2$
on $\A$ is a (group) homomorphism from $\bbZ^2$ to $Aut(\A).$ The triple $(\A, \bbZ^2, \ga)$ is called a $\cstar$-
dynamical system. In the case of $\bbZ^2,\ \ga$ defines two commuting automorphisms of $\A,\
\ga_1 = \ga(1, 0),\ \ga_2 = \ga(0, 1).$ Conversely, any two commuting automorphisms of $\A$ define a $\bbZ^2$-
action $\ga.$

We briefly recall the crossed product and semicrossed product constructions. (E.g.,
see \cite{pedersen},
 \cite{peters} for more details.) Let

\[ \il^1(\bbZ^2, \A, \ga) = \{ F:\ \bbZ^2 \to \A, \text{ such that } \sum_{i, j \in \bbZ} ||F(i,j)|| < \infty \} , \]

equipped with the usual multiplication

\[ (\gd_{ij}\otimes f)(\gd_{kl} \otimes g) = \gd_{i+k, j+l}\otimes f\ga(i,j)(g) = \gd_{i+k, j+l} \otimes f\ga_1^i
\ga_2^j(g) \]

and involution

\[ (\gd_{ij} \otimes f)^* = \gd_{-i, -j} \otimes f^* \]

where $f, g \in \A,$ and $\gd_{ij} \otimes f$ denotes the function $F: \bbZ^2 \to \A$ which takes the value $f$ at
the point $(i,j) \in \bbZ^2$ and zero elsewhere.

The crossed product, $\bbZ^2 \times_{\ga} \A,$ is the completion of $\il^1(\bbZ^2, \A, \ga)$ with respect to the norm
\begin{multline*}
||F|| = \sup \{ ||\pi(F)||: \pi \text{ a nondegenerate} \\
 \text{star representation of } \il^1(\bbZ^2, \A, \ga) \}.
\end{multline*}

The subalgebra $\il^1(\bbZ_+^2, \A, \ga)$ of $\il^1(\bbZ^2, \A, \ga)$ consists of
functions $F: \bbZ^2 \to \A$ such that $F(i,j) = 0 $ if either $i$ or $j$ is negative.
The subalgebra is not star-closed. The semicrossed product, $\scp \A,$ has been
defined in other contexts (e.g., \cite{peters}) as the completion of the
$\il^1$-algebra with respect to the norm
\[ ||F||_c := \sup \{ ||\pi(F)||: \pi \text{ a contractive representation of } \il^1(\bbZ_+^2, \A, \ga)\}.\]
Since nondegenerate star representations are contractive, we have $||F|| \leq ||F||_c.$ But in this case the
reverse inequality also holds.

\begin{lem} \label{l:scpnorm}
For any $F \in \il^1(\bbZ_+^2, \A, \ga),$ we have $||F||_c = ||F||,$ and consequently $\scp \A$ is the completion
of $\il^1(\bbZ_+^2, \A, \ga)$ in $\bbZ^2 \times_{\ga} \A.$
\end{lem}

\begin{proof}
Let $\rho$ be a contractive representation of $\il^1(\bbZ_+^2, \A, \ga).$ It is then
also contractive with respect to $||\cdot ||_c,$ and hence it can be extended to a
contractive representation of $\scp \A.$ By Ling and Muhly \cite{lingmuhly} $\rho$ is
completely contractive, and by Arveson-Stinespring theorem \cite{arveson}, there
exists a $\cstar$ dilation $\td{\rho}.$ If $\rho$ acts on the Hilbert space $\H,$ and
$\td{\rho}$ acts on $\K,$ there is an isometry $V$ from $\H$ into $\K$ such that
$\rho(F) = V^*\td{\rho}(F)V.$ Hence, $||\rho(F)|| \leq ||\td{\rho}(F)|| \leq ||F||.$

It follows that $||F||_c \leq ||F||.$
\end{proof}

\subsection{Canonical representations of semicrossed products} \label{ss:canrep}
By Pedersen (\cite{pedersen}), Theorem 7.7.7, and the fact that $\bbZ^2$ is amenable,
 it follows that the full crossed product
coincides with the reduced crossed product.  Hence if $(\pi, \H)$ is a faithful representation of $\A,$ then 
\linebreak[2]
\mbox{$(\td{\pi} \times U, \ell^2(\bbZ^2, \H))$} is a faithful representation of $\bbZ \times_{\ga}\A,$ where
$U : \bbZ^2 \to \L(\ell^2(\bbZ^2, \H)))$ \\
$U(s, t):= U_1^s U_2^t,$ and $U_1(\gd_{m,n}\otimes h) = \gd_{m+1,n}\otimes h,$ and $U_2(\gd_{m,n}\otimes h) = \gd_{m,n+1}\otimes
h.$ We think of $U_1,\ U_2$ as the horizontal and vertical shifts.

$\td{\pi}$ is a representation of $\A$ on $\ell^2(\bbZ^2, \H)$ given by
\[ \td{\pi}(a)(\gd_{m,n}\otimes h) = \gd_{m,n}\otimes \pi(\ga(-m, -n)(a) =
\gd_{m,n} \otimes \pi(\ga_1^{-m}\ga_2^{-n}a) .\]

Finally,  $\td{\pi}\times U$ is a representation of  $\bbZ_+^2 \times_{\ga} \A$ on $ \ell^2(\bbZ^2, \H)$ defined on
generating elements $\gd_{i,j}\otimes a$ by
\[ (\td{\pi}\times U)(\gd_{i,j}\otimes a) = \td{\pi}(a)U_1^iU_2^j .\]

In this paper we are concerned with actions of $\bbZ^2$ on finite dimensional $\cstar$-algebras, and we will identify
the resulting semicrossed products with matrix algebras of analytic functions.  We begin with the very simplest case,
$\A = \bbC.$

\begin{eg} \label{e:bidisk}
If the $\cstar$-algebra $\A = \bbC,$ the actions $\ga_1, \ga_2$ are trivial.  Take $\pi$ to be the identity representation
of $\bbC$ on (the Hilbert space) $\H = \bbC.$  Then $\td{\pi}(a)(\gd_{m,n} \otimes h) = \gd_{m,n} \otimes \pi(a)h.$
Let $W: L^2(\bbT^2) \to \ell^2(\bbZ^2)$ be the two-dimensional Fourier transform,
\[ W(\sum_{m,n \in \bbZ} x_{mn}z^m w^n ) = \sum_{m,n \in \bbZ} \gd_{m,n}\otimes x_{mn} . \]

The representation $\Ad{W}\circ (\td{\pi}\times U)$ takes the element \mbox{$\sum_{i, j = 0}^N \gd_{ij}\otimes a_{ij}$}
in the semicrossed product to the multiplication operator $M_p$ on $L^2(\bbT^2),$\linebreak[2]$ M_pg = pg,$ where
$p(z, w) = \sum_{i,j = 0}^N a_{mn}z^mw^n, \ (N \in \bbZ^+).$ It follows that the semicrossed product
$\scp \bbC,$ which is the closure of the set of elements of the above form, is identified with the closure of the
polynomials in two variables in the supremum norm, that is, the bidisk algebra, $\iAdd.$
Of course the crossed product is the $\cstar$-envelope, $C(\bbT^2).$
\end{eg}
\vspace{.15in}
Next we characterize the general case of semicrossed products \linebreak[2] $\scp C(X)$, where $X$ is finite.
Fix an integer $k  > 1$ and let $\gs_1, \gs_2$ be two commuting permutations of $\{0, 1, \dots, k-1 \}$. Let
$X = \{x_0, \dots, x_{k-1} \}$ and $\ga = (\ga_1, \ga_2)$ the $\bbZ^2$-action on $\A = C(X)$ given by
$\ga_i(f)(x_j) = f(x_{\gs_i(j)}),\ i = 1, 2$.

\begin{prop} \label{p:z2}
The algebra $\scp C(X)$ is unitarily equivalent to the operator algebra on $\oplus_{j=0}^{k-1} L^2(\bbT^2)$ with the
generators $V_1 = zP_1,\ V_2 = wP_2,$ and $D(f),\ f \in C(X),$ where $P_i$ are the permutation matrices
\[ P_i = \sum_{j=0}^{k-1} e_{j, \gs_i(j)} \ i = 1, 2 \text{ and } D(f) = \sum_{j=0}^{k-1} f(x_j) e_{j, j}, \
f \in C(X) .\]
\end{prop}

\begin{proof}
Let $\pi$ be the faithful representation of $C(X)$ on $\H = \bbC^k, \ $
\mbox{$\pi(f)h = (f(x_0)h_0, \dots, f(x_{k-1})h_{k-1}),$}
for $h = (h_0, \dots, h_{k-1}) \in \H, \ f \in C(X).$ Then $\td{\pi}$ is a representation of $C(X)$ on
$\td{\H} = \ell^2(\bbZ^2, \H),$ and
\[ \td{\pi}(f)(\gd_{m,n} \otimes h) = \gd_{m,n} \otimes \pi(\ga_1^m \ga_2^n f)h .\]

Let $\mathbf{h}_j = (0, \dots, 0, 1, 0 \dots, 0), \ j = 0, \dots k-1$ be the standard basis vectors for $\H$, and
let $W: \oplus_{j=0}^{k-1} L^2(\bbT^2) \to \ell^2(\bbZ^2, \H) $ be the unitary operator given by
\[ W(0, \dots, 0, z^mw^n, 0, \dots, 0) = \gd_{m,n} \otimes \mathbf{h}_{j'}\]
where if $z^mw^n$ occurs in the $j^{th}$ column on the left (with the first column
corresponding to $j = 0$) , then $j' = \gs_1^{-m}\gs_2^{-n}(j)$.

We now calculate
\begin{align*}
W^*U_1W(0, \dots, 0, z^m w^n, 0, \dots, 0) &= W^*U_1(\gd_{m,n} \otimes \mathbf{h}_{j'}) \\
                                           &= W^* (\gd_{m+1,n} \otimes \mathbf{h}_{j'}) \\
                                           &= (0, \dots, 0, z^{m+1} w^n, 0, \dots, 0)
\end{align*}
where the term $z^{m+1} w^n$ appears in the $\gs_1(j)$ column. In other words, \\
$W^*U_1W = z(\sum_{j=0}^{k-1} e_{j,\ \gs_1(j)}) = zP_1 = V_1$.\\
Similarly, $W^*U_2W = V_2$.

Moreover, for $f \in C(X),$ and $g_j \in L^2(\bbT^2),\ j = 0, \dots, k-1,$ it is straightforward that
\[ W^* \td{\pi}(f)W(g_0,\dots, g_{k-1}) = (f(x_0)g_0, \dots, f(x_{k-1})g_{k-1}),\]
so that $W^* \td{\pi}(f)W = D(f).$

\end{proof}

\begin{eg} \label{e:bk2}
Let $\gs_1 = \gs_2 =\gs = (k-1\ k-2\ \dots \ 1\ 0)$ be the forward shift on $\{0, 1, \dots, k-1\}$. Then
\begin{equation*}
P_1 = P_2 = \begin{pmatrix}
                          0 & 0 & \dots & 1 \\
                          1 & 0 & \dots & 0 \\
                          0 & 1 & \dots & 0 \\
                         \vdots & \vdots &  & \vdots \\
                          0 & 0 & 1 & 0
             \end{pmatrix}
\end{equation*}

It is straightforward to verify that the algebra generated by $V_1 = zP_1,\ V_2 = wP_2,$ and the set of all diagonal
matrices generate the subalgebra of $\B_{k,2}$ for which the functions $f_{ij}$ are polynomials in $z, w$.  Thus if
$\ga_1, \ga_2$ are the automorphisms of $C(X)$
corresponding to $\gs_1, \gs_2,$ it follows from Proposition \ref{p:z2} that the
semicrossed product $\scp C(X)$ is isometrically isomorphic to $\B_{k,2}.$
\end{eg}

\begin{cor}
Let $\A = M_n \otimes C(X)$ and $\ga = (\ga_1, \ga_2)$ be as in the previous example, so
$\ga_i(B\otimes f) = B\otimes \ga_i(f),\ i = 1, 2.$ Then $\scp \A$ is isometrically isomorphic to $M_n \otimes \B_{k,2}$.
\end{cor}
\begin{proof}
Clear.
\end{proof}

\begin{rem} \label{r:compisom}
In Example \ref{e:bk2} and in the Corollary the condition `isometrically isomorphic' can be replaced by
`completely isometrically isomorphic'. Indeed, this is immediate from the unitary equivalence in Proposition
\ref{p:z2}.
\end{rem}

\subsection{Tensor products of $\B_k$'s}
As the operator algebras $\B_k$ were identified as semicrossed products by
(transitive) $\bbZ$-actions on $\bbC^k$ (cf section \ref{s:matrixfunctions}), it may seem intuitively clear that the
tensor product of $\B_k $ with $ \B_{\ell}$
would correspond to semicrossed products by $\bbZ^2$ `product actions' on $\bbC^k \otimes \bbC^{\ell}$. We will in fact
establish such a result. First, however, we examine the maximal ideal structure of these tensor products.

Realizing the operator  algebra $\B_k$ (cf section \ref{s:matrixfunctions}) as a
$k\times k$ matrix function algebra acting on the direct sum
$\oplus_{j=0}^{k-1}L^2(\bbT)$, the algebraic tensor product $\B_k \otimes \B_{\ell}$
acts naturally on $\oplus_{j=0}^{k\ell -1} L^2(\bbT)$. The spatial tensor product
$\B_k \wh{\otimes} \B_{\ell}$ is the completion of the algebraic tensor product in
this representation.  Alternatively, $\B_k \wh{\otimes} \B_{\ell}$ can be viewed as
the completion of the algebraic tensor product in $M_{k\ell}\otimes \iAdd$.

As before, let $e_{ij},\ 0 \leq i, j \leq k-1$ denote matrix units for $M_k$, and similarly let
$\ol{e}_{i'j'},\ 0 \leq i', j' \leq \ell - 1$ denote matrix units for $M_{\ell}$. Thus, the algebraic tensor product
$\B_k \otimes \B_{\ell}$ is  spanned by linear combinations of elements of the form
\begin{equation*}
  e_{ij}\otimes \ol{e}_{i'j'}\otimes z^{i-j [\text{mod }k] +mk}w^{i'-j'[\text{mod }\ell] +n\ell},
 \quad m,\ n \in \bbZ,\ m,\ n \geq 0
\end{equation*}
where as before $i-j [\text{mod } k]$ lies in $\{0, \dots, k-1\}$ \ (resp. $i'-j'
[\text{mod }\ell] \in \{0,\dots, \ell -1 \}$. Thus, $\B_k \wh{\otimes} \B_{\ell}$ is a
$k\ell \times k\ell$ matrix of functions $f_{(ij)(i'j')} \in \iAdd $ whose Fourier
series have the form
\[ f_{(ij)(i'j')}(z, w) \sim \sum_{m, n \geq 0}a_{mn} z^{i-j [\text{mod }k] +mk}w^{i'-j'[\text{mod }\ell] +n\ell}\].

Observe that functions in the diagonal subalgebra of $\B_k \wh{\otimes} \B_{\ell}$ belong to
$\iAkd \wh{\otimes}\mathit{A}_{\ell}(\mathbb{D})$.
\begin{lem}
The maximal ideals of $\iAkd \wh{\otimes}\mathit{A}_{\ell}(\mathbb{D})$ are kernels of the evaluation
homomorphisms $f \mapsto f(\gl, \mu),\ (\gl, \mu) \in \ol{\bbD}^2$.
\end{lem}
\begin{proof} This is in the spirit of Lemma \ref{l:iAkddmaxid}.
\end{proof}

As with the earlier Lemma, the correspondence between maximal ideals and points in the closed bidisk is not bijective.

For $0 \leq i_0 \leq k-1$ and $0 \leq i_0' \leq \ell -1$, define $\gp_{00}^{i_0i_0'}: \B_k \wh{\otimes} \B_{\ell}
\to \bbC$ by $\gp_{00}^{i_0i_0'}(F) = f_{(i_0i_0)(i_0'i_0')}(0,0)$, where $F = (f_{(ij)(i'j')})$. Define
$\gp_{0\mu}^{i_0}: \B_k \wh{\otimes} \B_{\ell} \to M_{\ell},
\ \gp_{0\mu}^{i_0}(F) = (f_{(i_0i_0)(i'j')}(0,\mu))_{i'j'=0}^{\ell -1},\quad
\gp_{\gl 0}^{i_0'}: \B_k \wh{\otimes} \B_{\ell} \to M_{k},\
\gp_{\gl 0}^{i_0'}(F) = (f_{(ij)(i_0'i_0')}(\gl,0))_{ij=0}^{k-1}, \quad \text{and }
\gp_{\gl \mu}: \B_k \wh{\otimes} \B_{\ell} \to M_{k\ell},\
\gp_{\gl \mu}(F) = F(\gl, \mu) = (f_{(ij)(i'j')}(\gl, \mu))_{(ij)(i'j')=(0)(0)}^{(k-1)(\ell -1)}.$

Denote the kernels of these homomorphisms, respectively, by \\
$J_{00}^{i_0i_0'}, \ J_{0\mu}^{i_0},\ J_{\gl 0}^{i_0'}$ and $ J_{\gl \mu}$.

\begin{prop} \label{p:maxid2}
Every maximal ideal of $\B_k \wh{\otimes} \B_{\ell}$ is of the form \\
$J_{00}^{i_0i_0'}, \ J_{0\mu}^{i_0},\
J_{\gl 0}^{i_0'}$ or $ J_{\gl \mu}$, for $\gl,\ \mu \in \ol{\bbD} \bsl \{0\}$.
\end{prop}

\begin{proof}  Note that the $\gp$-homomorphisms are
all onto simple algebras, so their kernels are maximal ideals.

The proof is similar to that of Proposition \ref{p:bk2maxid}. Let $J$ be a maximal ideal of $\B_k \wh{\otimes} \B_{\ell}$.
 Since $J$ cannot contain the identity, there are indices
$i_0, i_0'$ such that for all $F \in J$, $F = (f_{(ij)(i'j')}),\ f_{(i_0i_0)(i_0'i_0')}) \neq 1.$ Indeed, if that were
not the case, so we could find $F_{ii'}$ whose $((ii), (i'i'))$ component was $1,\ 0 \leq i \leq k-1, \quad
0 \leq i' \leq \ell -1$, then we would obtain
\[ \sum_{ i, i' = 0}^{i=k-1, i'= \ell-1} (e_{ii}\otimes \ol{e}_{i'i'})F_{ii'}(e_{ii}\otimes \ol{e}_{i'i'}) = I \in J \]
a contradiction.

Without loss of generality, we may assume that for all $F \in J,\ F = (f_{(ij)(i'j')}),\ f_{(00)(00)} \neq 1.$
Let $\I$ be the ideal in $\iAkd \wh{\otimes}\mathit{A}_{\ell}(\mathbb{D}) $ generated by
$\{ f = f_{(00)(00)}, \text{ for some } F = (f_{(ij)(i'j'))} \in J \}.$ By assumption $\I$ is proper, so it is
contained in some maximal ideal given by evaluation at $(\gl, \mu) \in \ol{\bbD}^2$.

If $(\gl, \mu) = (0, 0)$ then $J \subset J_{00}^{00}$ and hence $J = J_{00}^{00}$ by maximality of $J$.
If $(\gl, \mu) = (\gl, 0)$ for some $\gl \in \ol{\bbD} \bsl \{0\}$, we claim $J = J_{\gl,0}^{0}$.
Let $F \in J,\ F = (f_{(ij)(i'j')}),$ and select $ 0 \leq i, j \leq k-1$. From the ideal property of $J$ we obtain
\[ (z^{k-i}\otimes e_{0i}\otimes \ol{e}_{00})F (z^j \otimes e_{j0} \otimes \ol{e}_{00})
= z^{k+j-i} f_{(ij)(00)} \otimes e_{00} \otimes \ol{e}_{00} \in J \]
whence $\gl^{k+j-i}f_{(ij)(00)}(\gl, 0) = 0$, and in particular $f_{(ij)(00)}(\gl, 0) = 0$. This shows that
$J \subset J_{\gl,0}^{0}$, and by maximality of $J$ equality prevails.

The other cases are analogous.

\end{proof}

\begin{cor} \label{c:notisom}
For positive integers $k, \ell, p,\ k, \ell > 1,$ the Banach algebras $\B_{p,2},\ \B_k \wh{\otimes} \B_{\ell}$
are not isomorphic.
\end{cor}

\begin{proof}
$\B_{p,2}$ has maximal ideals of at most two distinct codimensions, whereas $ \B_k \wh{\otimes} \B_{\ell}$
has maximal ideals of at least three distinct codimensions.
\end{proof}

\subsection{Perpendicular actions} \label{ss:perpact}
In example \ref{e:bk2} the algebra $\B_{k,2}$ was identified with the semicrossed product $\scp C(X)$, where
$\ga = (\ga_1, \ga_2)$ in which $\ga_1 = \ga_2$ was implemented by a cyclic permutation acting transitively on
$X$. In this subsection we consider the opposite extreme: namely, where $\ga = (\ga_1, \ga_2),\ \ga_i$ is implemented
by a permutation $\gs_i,\ i = 1, 2$ such that $\gs = (\gs_1, \gs_2)$ is transitive as a $\bbZ^2$ action, but neither
$\gs_i$ is transitive, and furthermore the orbits of $\gs_1,\ \gs_2$ overlap minimally, that is, in a single point.
To state this formally, we make the following

\begin{defn} \label{d:perpact}
Let $X$ be a finite set, $\gs_i, \ i = 1, 2$ commuting permutations on $X$, and $\gs = (\gs_1, \gs_2)$ the
induced $\bbZ^2$-action on $X$.  We say that $\gs$ is a \emph{perpendicular} action if \\
There exist finite sets $X_1,\ X_2$, cyclic (transitive) permutations $\gt_i$ on $X_i \ (i = 1, 2)$ and
a bijection $h: X \to X_1 \times X_2$ such that for all $x_0 \in X,$
\begin{multline*}
 h(\gs_1(x_0)) = (\gt_1(x_1), x_2), \text{ and } h(\gs_2(x_0)) = (x_1, \gt_2(x_2)), \\
 \text{where } (x_1, x_2) = h(x_0) .
\end{multline*}
In other words, $h\circ \gs = \gt \circ h,$ where $\gt(m, n) = \gt_1^m \times \gt_2^n$ on $X_1 \times X_2$.
\end{defn}

In the language of dynamics, $\gs$ is a perpendicular action if it is conjugate to a product of transitive actions
on finite sets.

While the following proposition is elementary, we are not aware of it in the dynamical systems literature.
\begin{prop} \label{p:peract}
Let $X$ be a finite set, $\gs_i \ (i = 1, 2)$ commuting permutations, and $\gs = (\gs_1, \gs_2)$ the
induced $\bbZ^2$-action.
Then $\gs$ is a perpendicular action if and only if
\begin{enumerate}
\item
 For any $x, y \in X,\ \O_{\gs_1}^x \cap \O_{\gs_2}^y$ consists of a single point, where $\O_{\gs_1}^x \
(\text{resp. } \O_{\gs_2}^y$) denotes the orbit of $x$ under $\gs_1$ (resp., the orbit of $y$ under $\gs_2$).

\item
All $\O_{\gs_1}^x \ (x \in  X)$ have the same cardinality, and all
$\O_{\gs_2}^y \ (y \in Y)$ have the same cardinality.
\end{enumerate}
\end{prop}

\begin{proof}
First, if $\gs$ is a perpendicular action it is transitive, and it is clear that (i) and (ii) are satisfied.

Suppose now that (i)  and (ii) are satisfied. Observe that $\gs$ is a transitive action.  Indeed, if
$x, y \in X$ by (i) there are $m, n \in \bbZ$ such that $\gs_1^m(x) = \gs_2^n(y);$ that is
$ y = \gs(m, -n)(x)$, which is transitivity.

Choose $x_0 \in X$ and let $X_1 = \O_{\gs_1}^{x_0}$ and $X_2 = \O_{\gs_2}^{x_0}$. Define a map
$h: X \to X_1 \times X_2$ by $h(x) = (y_1, y_2)$ where $y_1$ is the unique element in
$\O_{\gs_1}^{x_0} \cap \O_{\gs_2}^x$, and $y_2$ is the unique element in
$\O_{\gs_1}^x \cap \O_{\gs_2}^{x_0}.$ We show that the mapping $h$ is a bijection.

Injectivity: Let $h(x) = h(x')$. Then $(y_1, y_2) = (y_1', y_2')$ where $\O_{\gs_1}^{x_0} \cap \O_{\gs_2}^x =
\{ y_1 \} = \{ y_1'\} = \O_{\gs_1}^{x_0} \cap \O_{\gs_2}^{x'}.$ Since the orbits $\O_{\gs_2}^x$ and
$\O_{\gs_2}^{x'}$ have a point in common, they are equal.  Similarly, $\O_{\gs_1}^x = \O_{\gs_1}^{x'}.$
So we have $\O_{\gs_1}^{x'} \cap \O_{\gs_2}^x =  \O_{\gs_1}^x \cap \O_{\gs_2}^x = \{ x\}$, and on the other hand
$\O_{\gs_1}^{x'} \cap \O_{\gs_2}^x = \O_{\gs_1}^{x'} \cap \O_{\gs_2}^{x'} = \{x' \}.$ Hence $ x = x'.$

Surjectivity: Let $(y_1, y_2) \in X_1 \times X_2.$ Let $z$ denote the unique element in the intersection
$\O_{\gs_1}^{y_2} \cap \O_{\gs_2}^{y_1}$. As $y_1$ is in $\O_{\gs_2}^z$ and $y_2 $ is in $\O_{\gs_1}^z,$ we have
$\O_{\gs_1}^{x_0} \cap \O_{\gs_2}^z = \{ y_1\}$ and $\O_{\gs_1}^z \cap \O_{\gs_2}^{x_0} = \{y_2 \}$. Hence
$h(z) = (y_1, y_2)$.

Let us now show that the action is perpendicular.  Take $(y_1, y_2) \in X_1 \times X_2 $ and set
$\wh{\gs}_i := h^{-1} \circ \gs_i \circ h,\ i = 1, 2.$ Let $ h = (h_1, h_2);$ that is,
$h_i = \pi_i \circ h,$ where $\pi_i$ is the $i^{th}$ coordinate projection, $i = 1, 2.$
First we show that $\pi_2 \circ \wh{\gs}_1(y_1, y_2) = y_2.$ Call $ z = h^{-1}(y_1, y_2)$.
As $\O_{\gs_1}^z = \O_{\gs_1}^{\gs_1(z)}$, it follows $ \{ h_2(\gs_1(z))\}
= \O_{\gs_1}^{\gs_1(z)} \cap \O_{\gs_2}^{x_0} = \O_{\gs_1}^z \cap \O_{\gs_2}^{x_0} = \{ h_2(z) \}$.
Hence, $h_2(\gs_1(z)) = h_2(z)$. Therefore there is a map $\gt_1: X_1 \times X_2 \to X_1$ such that
$\wh{\gs}_1(y_1, y_2) = (\gt_1(y_1, y_2), y_2)$.

It remains to show that $\gt_1$ depends only on $y_1$, and hence can be viewed as a map $X_1 \to X_1$.
We show that $\gt_1(y_1, y_2) = \gt_1(y_1, y_2')$. Let $h^{-1}(y_1, y_2) =z,\ h^{-1}(y_1, y_2') = z'$.
Since $\gs$ is a transitive action, $z' = \gs(m,n)z = \gs_1^m \gs_2^n(z)$ for some $(m, n) \in \bbZ^2$.
But as we have noted above, $h_2(z) = h_2(\gs_1(z))$ and similarly $h_1(z) = h_1(\gs_2(z))$,
 so it is an action of $\gs_2$ that maps $z$ to $z'$,
say $z' = \gs_2^n(z)$. Since $\gs_1,\ \gs_2$ commute, $\gs_1(z') = \gs_1(\gs_2^n(z)) = \gs_2^n(\gs_1(z))$,
so that $\gs_1(z'),\ \gs_1(z)$ are on the same $\gs_2$-orbit. Hence, $\gt_1(y_1, y_2) = h_1(\gs_1(z')) =
\gt_1(y_1, y_2')$.

\end{proof}

\begin{rem} \label{r:perpact}
1. If $\gs$ is a perpendicular $\bbZ^2 $ action which is conjugate to a product action
$\gt_1 \times \gt_2$ on $X_1 \times X_2$ and also conjugate to $\gt_1' \times \gt_2'$ on $X_1' \times X_2'$,
then card$(X_i)$ = card$(X_i'),\ i = 1, 2$. \\
2. Is there a version of Proposition \ref{p:peract} for topological dynamics?
\end{rem}

\begin{eg} \label{e:perpact}
Take $X = \{0, 1 \dots, 5\},\ \gs_1 = (0\, 1\, 2)(3\, 4\, 5)$ and $ \gs_2 = (0\, 3)(1\, 4)(2\, 5)$.
Then the conditions of the Proposition are satisfied, and with $X_1 = \{0, 1, 2\},
X_2 = \{0, 3\}, $ and $ x_0 = 0$ we have
$h(i) = (i, 0),\ i = 0, 1, 2$ and $ h(i) = (i, 3), \ i = 3, 4, 5 $.
\end{eg}

\subsection{Semicrossed Products by perpendicular actions} \label{ss:sppa}
Earlier (in \ref{ss:canrep}) we reviewed the construction of the (semi)crossed product as the completion of
an $\ell_1$-algebra in a canonical representation. In the case of a perpendicular $\bbZ^2$-action, the
canonical representation can be taken to be a product representation, and this allows us to view the
(semi)crossed product as the tensor product of (semi)crossed products with respect to $\bbZ$-actions. In what
follows we consider (semi)crossed products of $C(X)$ where $X$ is a compact metric space.

Let $X_i$ be a compact metric space, $  \gs_i$ a homeomorphism of $X_i$, and $\ga_i$ the induced automorphism
of $C(X_i),\ \ga_i(f) = f\circ \gs_i,\ f \in C(X_i),\ i =1, 2.$ Let $\pi_i$ be a faithful representation of
$C(X_i)$ on a Hilbert space $\H_i$. Denote by $K(\bbZ, C(X_i), \ga_i)$ the dense
subalgebra of $\ell_1(\bbZ, C(X_i), \ga_i)$
consisting of finite sums of the generators $\sum_m \gd^{i}_m \otimes f_m$. Then, as in \ref{ss:canrep},
the (semi)crossed product $\bbZ \times_{\ga_i} C(X_i)$ (resp., $\bbZ^+\times_{\ga_i} C(X_i))$
is the completion of $K(\bbZ, C(X_i), \ga_i)$ (resp., $K(\bbZ^+, C(X_i), \ga_i)$)
in the representation $\td{\pi_i} \times U_i$ on the Hilbert space $\td{\H}_i = \ell^2(\bbZ, \H_i)$.

\begin{multline} \label{E:tensor}
(\td{\pi}_1 \times U_1)\otimes (\td{\pi}_2 \times U_2)(K(\bbZ, C(X_1), \ga_1))\otimes (K(\bbZ, C(X_2), \ga_2))\\
= (\td{\pi}_1\otimes \td{\pi}_2)\times U)(K(\bbZ^2, C(X_1)\otimes C(X_2), \ga))
\end{multline}
where $\ga$ is the $\bbZ^2$-action $\ga = (\ga_1, \ga_2), \ \ \td{\pi}_1\otimes \td{\pi}_2$ acts of the
Hilbert space $\td{\H} = \td{\H}_1 \wh{\otimes} \td{\H}_2$ and $U(m,n) = U_1^m \otimes U_2^n $.
Since the algebraic tensor product $C(X_1) \otimes C(X_2)$ is dense in $C(X_1 \times X_2)$, and
$\pi_1 \otimes \pi_2$ is a faithful representation of  $C(X_1 \times X_2)$, it follows that the completion
of the right-hand side of (\ref{E:tensor}) is the crossed product $\cp C(X), \ X = X_1 \times X_2$.

\begin{prop} \label{p:prodact}
With notation as above,
\[ (\bbZ \times_{\ga_1} C(X_1)) \wh{\otimes} (\bbZ \times_{\ga_2} C(X_2)) = \cp C(X) \]
and
\[ (\bbZ_{+} \times_{\ga_1} C(X_1)) \wh{\otimes} (\bbZ_{+} \times_{\ga_2} C(X_2)) = \scp C(X) .\]
\end{prop}

\begin{proof} The completion of $(K(\bbZ, C(X_1), \ga_1))\otimes (K(\bbZ, C(X_2), \ga_2))$ in the norm provided
by the representation $(\td{\pi}_1 \times U_1)\otimes (\td{\pi}_2 \times U_2)$ on $\L(\td{\H})$ is a C$^*$-norm on
the algebraic tensor product $(\bbZ \times_{\ga_1} C(X_1)) \otimes (\bbZ \times_{\ga_2} C(X_2))$.
 However, the crossed product $\bbZ \times_{\ga_i} C(X_i)$
is nuclear (e.g., \cite{rordam}, Proposition 2.1.2), so there is only one completion, and hence
the completion is $\cp C(X)$. The statement about semicrossed products follows from the fact that the
semicrossed product norm is, by Lemma \ref{l:scpnorm}, the restriction of the C$^*$-crossed product norm.

\end{proof}

\begin{cor} \label{c:perpact}
Let $X$ be a finite set, $\gs$ a perpendicular $\bbZ^2$-action on $X$, so that $\gs$ is conjugate
to $\gt_1 \times \gt_2$ acting on $X_1 \times X_2$, where $\gt_i$ is a transitive $\bbZ$-action on $X_i$.
If card$(X_1) = k,\ $ card$(X_2) = \ell$, then the semicrossed product
$\bbZ_{+}^2 \times_{\gs} C(X)$ is identified with $\B_k \wh{\otimes} \B_{\ell}$.

\end{cor}

\begin{proof} This follows immediately from the above and the identification of $\bbZ_{+} \times_{\gs_i} C(X_i)$
with $\B_k \  (i=1)$ or $\B_{\ell} \ (i=2)$.
\end{proof}

\begin{cor} With the same assumptions as in Corollary \ref{c:perpact}, $\scp C(X)$
is  generated by $\{zP_k \otimes wP_{\ell}, D_k \otimes D_{\ell} :\ D_k \in \text{diag}(M_k),
D_{\ell} \in \text{diag}(M_{\ell}) \}$, where $\text{diag}(M_k)$ (resp., $\text{diag}(M_{\ell})$) is the
algebra of diagonal matrices in $M_k$ (resp., $M_{\ell}$.)

$P_k \in M_k,$ and $P_{\ell} \in M_{\ell}$ have the form of the matrix in Example \ref{e:bk2}
\end{cor}

\begin{proof}
This follows from the form of $\B_k$ (See (\ref{E:Ak})) and the fact that the subalgebra of $\B_k$
consisting of matrix functions with polynomial entries are dense in $\B_k$.
\end{proof}

\begin{cor}
There are transitive $\bbZ^2$-actions $\gs,\ \gs'$ on a finite set $X$ such that the semicrossed products
$\bbZ^2 \times_{\gs} C(X)$ and $\bbZ^2 \times_{\gs'} C(X)$ are not isomorphic.
\end{cor}
\begin{proof}
This follows from Corollary \ref{c:notisom}, Example \ref{e:bk2}, and Corollary \ref{c:perpact}.
\end{proof}

Acknowlegements: \ The first author would like to thank \ Prof. Peters for
his invitation to Iowa State university where part of this work has been
done . He would like also to acknowledge the support of KFUPM\ where the
second part has been finished.

\end{document}